\documentclass[11pt]{amsart}
\usepackage{amsmath}
\usepackage{amsxtra}
\usepackage{amscd}
\usepackage{amsthm}
\usepackage{amsfonts}
\usepackage{amssymb}
\usepackage{eucal}

\numberwithin{equation}{subsection}

\begin{document}

\newcommand{\thmref}[1]{Theorem~\ref{#1}}
\newcommand{\secref}[1]{\S~\ref{#1}}
\newcommand{\defref}[1]{Definition~\ref{#1}}
\newcommand{\lemref}[1]{Lemma~\ref{#1}}
\newcommand{\propref}[1]{Proposition~\ref{#1}}
\newcommand{\corref}[1]{Corollary~\ref{#1}}
\newcommand{\remref}[1]{Remark~\ref{#1}}

\newcommand{\nc}{\newcommand}

\nc{\mc}{\mathcal}
\nc{\on}{\operatorname}
\nc{\Z}{{\mbb Z}}
\nc{\C}{{\mbb C}}
\nc{\Oo}{{\mc O}}
\nc{\D}{{\mc D}}
\nc{\E}{{\mc E}}
\nc{\Ll}{{\mc L}}
\nc{\M}{{\mc M}}
\nc{\la}{\lambda}
\nc{\ep}{\epsilon}
\nc{\V}{\mc{V}}
\nc{\inv}{^{-1}}
\nc{\ol}{\overline}
\nc{\wt}{\widetilde}
\nc{\wh}{\widehat}
\nc{\mb}{\mathbf}
\nc{\mf}{\mathfrak}
\nc{\mbb}{\mathbb}
\nc{\Cx}{{\mbb C}^\times}
\nc{\un}{\underline}
\nc{\sm}{\setminus}
\nc{\bs}{\backslash}
\nc{\Mg}{{{\mf M}_g}}
\nc{\Proj}{{\mc P}roj}
\nc{\Conn}{{\mc C}onn}
\nc{\excon}{{\mc E}x\Conn}
\nc{\exhiggs}{{\mc E}x{\mc H}iggs}
\nc{\boxt}{\boxtimes}
\nc{\ot}{\otimes}
\nc{\al}{\alpha}
\nc{\pa}{\partial}
\nc{\T}{{\mc T}}
\nc{\Ohalf}{\Omega_X^{\frac{1}{2}}}
\nc{\Sf}{{\mbb S}}
\nc{\vphi}{\varphi}
\nc{\io}{\iota}
\nc{\Mf}{{\mf M}}
\nc{\Cc}{{\mc C}}
\nc{\Pp}{{\mbb P}}
\nc{\Ev}{E^{\vee}}
\nc{\Mo}{{\mc M}^{\circ}}
\nc{\Mobar}{\ol{{\mc M}}^{\circ}}
\nc{\s}{{\mf s}}
\nc{\sbar}{{\ol{\s}}}
\nc{\Ac}{{\mc A}}
\nc{\XX}{X\times X}
\nc{\K}{\mf K}
\nc{\Kl}{\mbb K}
\nc{\W}{\mbb W}
\nc{\Ptil}{{\mc P}^{g-1}}
\nc{\Pc}{{\mc P}}
\nc{\Pg}{\on{Pic}_X^{g-1}}
\nc{\Jac}{\on{Jac}_X}
\nc{\Ho}{{\rm H}^0}
\nc{\Hi}{{\rm H}^1}
\nc{\Kern}{{\mc Kern}}
\nc{\Ebar}{\ol{E}}
\nc{\Oxy}{{\Oo({\mf y}-{\mf x})}}
\nc{\arrtheta}{\stackrel{\rightarrow}{\theta}}
\nc{\Op}{{\mc Op}}
\nc{\Mop}{{\mc M}{\mc Op}}
\nc{\bra}{\langle}
\nc{\ket}{\rangle}
\nc{\N}{{\mf N}}
\nc{\sw}{{\mf sl}}
\nc{\gl}{{\mf gl}}
\nc{\Sug}{{\Bbb S}}

\title{Opers and Theta Functions}

\author{David Ben-Zvi}

\address{Department of Mathematics, The University of Chicago,
5734 University Avenue, Chicago, Illinois 60637, USA}

\email{benzvi@math.uchicago.edu}

\author{Indranil Biswas}

\address{School of Mathematics, Tata Institute of Fundamental
Research, Homi Bhabha Road, Bombay 400005, INDIA}

\email{indranil@math.tifr.res.in}

\date{}

\begin{abstract}
We construct natural maps (the Klein and Wirtinger maps) from moduli
spaces of vector bundles on an algebraic curve $X$ to affine spaces,
as quotients of the nonabelian theta linear series. We prove a
finiteness result for these maps over generalized Kummer varieties
(moduli of torus bundles), leading us to conjecture that the maps are
finite in general.  The conjecture provides canonical explicit
coordinates on the moduli space. The finiteness results give
low--dimensional parametrizations of Jacobians (in $\Pp^{3g-3}$ for
generic curves), described by $2\Theta$ functions or second
logarithmic derivatives of theta.

We interpret the Klein and Wirtinger maps in terms of opers on
$X$.  Opers are generalizations of projective structures, and can be
considered as differential operators, kernel functions or special
bundles with connection. The matrix opers (analogues of opers for
matrix differential operators) combine the structures of flat vector
bundle and projective connection, and map to opers via generalized
Hitchin maps.  For vector bundles off the theta divisor, the Szeg\"o
kernel gives a natural construction of matrix oper. The Wirtinger map
from bundles off the theta divisor to the affine space of opers is
then defined as the determinant of the Szeg\"o kernel. This
generalizes the Wirtinger projective connections associated to theta
characteristics, and the assoicated Klein bidifferentials.

\end{abstract}

\maketitle

\section{Introduction.}

Let $X$ be a compact connected Riemann surface (or equivalently, a
connected smooth projective algebraic curve over $\C$).  Let
$\Mf_X(n)$ denote the moduli space of semistable vector bundles over
$X$ of rank $n$ and Euler characteristic zero (hence of degree
$n(g-1)$), and $\N_X(n)\subset\Mf_X(n)$ is the moduli space of vector
bundles with fixed determinant $\Omega_X^{\frac{n}{2}}$ (for a fixed
theta--characteristic $\Ohalf$ on $X$).  Let $\Theta\subset \M_X(n)$
denote the canonical theta divisor; its complement $\M_X(n)\sm\Theta$
is an affine variety, parametrizing rank $n$ vector bundles with
vanishing cohomology.

The theory of nonabelian theta functions provides an embedding of the
$(n^2-1)(g-1)$--dimensional affine variety $\N_X(n)\sm\Theta$ into an
affine space of dimension $n^g$. Specifically, by restricting the
canonical theta function to the image of the Jacobian $\Jac$ in
$\Mf_X(n)$ obtained by translating a given $E\in \N_X(n)\sm\Theta$ by
line bundles, we obtain elements of the $n\Theta$ linear series on the
Jacobian.

It is tempting to look for lower--dimensional parametrizations of
$\N_X(n)\sm\Theta$ which come closer to giving explicit coordinates on
the moduli space. Optimistically, one can hope for a natural finite
map of $\N_X(n)\sm \Theta$ to affine space of the same dimension. In
this paper we give a construction of such a map to affine space, which
we conjecture is finite, and explain its relations to theta functions,
projective structures, and differential operators on the Riemann
surface.

Our construction assigns special differential operators, or {\em
opers}, on $X$ to a vector bundle $E$ with vanishing
cohomologies. We define a map, the {\em Wirtinger map}
$\W$, from $\N_X(n)\sm\Theta$ to the space of all opers, which is an
affine space for the Hitchin base space of $X$. The dimension of
the space of opers is same as that of the moduli space (namely,
$(n^2-1)(g-1)$). By realizing the opers as {\em kernel functions}
on $X\times X$ we define the {\em Klein map} $\Kl$, 
sending the moduli space to a
(somewhat bigger) affine space of bidifferentials. Our main result
establishes the finiteness of the Klein map (for all $X$) and the
Wirtinger map (for generic $X$) when restricted to the moduli space of
torus bundles -- the generalized Kummer variety
$\K_X(n)\subset\N_X(n)$.

The case $n=2$ provides new finite parametrizations of $\Jac\sm\Theta$
(factoring through the Kummer $\K_X(2)=\Jac/\{\Ll\sim\Ll^*\}$) in
affine spaces of dimensions $\binom{g}{2}$ and $3g-3$ (that is,
quadratic and linear in the genus $g$), improving on the
parametrization given by the $2\Theta$ linear series (which requires
{\em exponential} dependence $2^g$ on the genus to embed the Kummer).
As a side--note we obtain that the collection of second logarithmic
derivatives of the theta function (considered in \cite{Tata I})
suffice to give a (generically) finite parametrization of the
Jacobian, and hence of a generic abelian variety.  Our proof uses the
behavior of the (abelian) Szeg\"o kernel near the theta divisor (in
fact near blowups of Brill--Noether loci) to show that the maps are
proper, hence finite, on the affine varieties $\K_X(n)\sm \Theta$
(giving finite extensions of the Gauss map of the theta divisor).

The Klein and Wirtinger maps may be defined either in terms of
restrictions of theta functions, or in terms of determinants of
nonabelian Szeg\"o kernels. From the point of view of theta functions,
the maps appear as certain quotients of the theta linear series,
obtained by restricting the theta function first from $\N_X(n)$ to
$\Jac$, then to $X\times X$ (via the difference map that sends $(x\,
,y)\in X\times X$ to ${\mathcal O}_X(x-y)$) and further to the $n$th
order infinitesimal neighborhood of the diagonal. The theta function
thereby defines kernel functions, sections on $X\times X$ of certain
sheaves of differentials. Such kernel functions, expanded near the
diagonal, are naturally interpreted as differential operators acting
between different line bundles on $X$. On the $n$th order
infinitesimal neighborhood of the diagonal, we obtain monic
differential operators with vanishing subprincipal symbol, which we
interpret as $\text{SL}_n$--opers on $X$.

Opers (for a reductive group
$G$) are special principal bundles with connection, which play a
central r\^ole in integrable systems and representation theory of loop
algebras. They were introduced in \cite{opers} in the context of the
geometric Langlands program, providing a coordinate--free expression
for the connections which appeared first in \cite{DS} as the phase
space of the generalized Korteweg--de Vries hierarchies. Opers form an
affine space, modeled on the vector space which is the base of
Hitchin's integrable system on the cotangent bundle of the moduli
space of bundles. For $G$ a classical group, opers are identified
with certain differential operators acting between line bundles on
$X$. In the case $G=\text{SL}_2$, opers are identified (after the
choice of a theta characteristic $\Ohalf$, which we fix) with
projective connections (or projective structures) on $X$.

By writing opers in terms of their kernel function, we obtain explicit
constructions of opers, generalizing the constructions of projective
structures from theta functions due to Klein and Wirtinger
(\cite{Tyurin}). This helps clarify some constructions of differential
operators on Riemann surfaces with projective structure (\cite{BR}).

Another point of view on the Klein and Wirtinger maps is given by
matrix opers and the Szeg\"o kernel. We define matrix opers by
applying the oper interpretation of differential operators to matrix
differential operators. Matrix opers combine the structures of
connections on a vector bundle and oper in a natural way (they play
the same r\^ole for multicomponent soliton equations that opers play
for KdV). A special class of matrix opers, the extended connections
(combining connections with projective structures) appear in \cite{BS}
(implicitly) and \cite{part 1} (explicitly) as twisted cotangent
spaces to the universal moduli space of vector bundles on Riemann
surfaces.  In analogy with the Hitchin system, we may apply invariant
polynomials to matrix opers, and obtain (scalar) opers. For extended
connections, we show the determinant map in fact gives a deformation
of the quadratic Hitchin map, which appears in the theory of
Virasoro--Kac--Moody algebras and isomonodromic deformation
(\cite{Sugawara}).

To every vector bundle $E$ off of the theta divisor, there is a
canonical matrix oper on $E$, defined by the nonabelian Szeg\"o kernel
of Fay (\cite{Fay nonabelian, Fay Szego, part 1}).
Applying the determinant map
to the Szeg\"o kernel we recover the pullback of the theta function,
and thus the Klein and Wirtinger maps. This point of view is motivated
by conformal field theory, where this map appears from taking
correlation functions associated to $\mc W$--algebra symmetries of
current algebras. We hope to describe this point of view in future
work, and expect it to facilitate the precise description of the
Wirtinger map and the proof of our finiteness conjecture.

Since we believe the point of view provided by opers is important in
understanding the r\^ole of Klein and Wirtinger maps, we describe
their structure in some detail in the first section. However, we
recommend readers to first jump ahead to the last two sections (which
can be read largely independently) where the maps are described in
elementary terms. The paper is organized as follows: in \secref{diffops as
kernels} we review the description of differential operators as kernel
functions, review some basics of opers, and describe matrix opers,
extended connections and their analogue of the Hitchin map. In
\secref{Klein and Wirtinger} we introduce the Klein and Wirtinger maps
via the Szeg\"o kernel, and prove the finiteness theorem for Kummers,
\thmref{the theorem}. Finally in \secref{thetas} we explain the
relation with classical constructions with theta functions, and draw
conclusions about $2\Theta$ functions and logarithmic derivatives on
Jacobians.

\section{Differential operators and kernels}\label{diffops as kernels}

\subsection{Notations.}
Let $X$ be a compact connected Riemann surface of genus
$g$ (a connected smooth complex projective
curve -- unless
explicitly noted, all constructions will be algebraic). Let $p_i\,
:\, X\times X\, \longrightarrow\, X$, $i=1,2$, be the projection to
the $i$-th factor. The diagonal divisor on $X\times X$ will be denoted
by $\Delta$. The involution on $X\times X$ given by interchange of
factors will be denoted by $\sigma$, so $\sigma(x,y)=(y,x)$. Given
holomorphic vector bundles $V$ and $W$ on $X$, we denote vector
bundles on $X\times X$ by
$$V\boxt W\, :=\, p_1^*V\otimes p_2^* W, \hskip.2in 
V\boxt W(n\Delta) \,:=\, p_1^* V\otimes p_2^* W\otimes\Oo_{X\times
X}(n\Delta)
$$ In particular $p_1^*V=V\boxt \Oo$ and $p_2^*W=\Oo\boxt W$. For a
vector bundle $V$ on $X$ we denote by $V^{\vee}=V^*\otimes \Omega_X$
the Serre dual vector bundle, where $\Omega_X$ is the holomorphic
cotangent bundle of $X$.  For a sheaf $W$, we will denote by $\Gamma
(W)=H^0(X,W)$ and $h^i(W)=\dim H^i(X, W)$.

Given a holomorphic vector
bundle $V$ over a complex manifold $M$, a \textit{torsor}, or
\textit{affine bundle}, for $V$ over
$M$ is a submersion of complex manifolds $\pi:A\to M$ with a simply
transitive, holomorphic action of the sheaf of sections of $V$ on the
sections of $A$. So the map $A\times_M V\, \longrightarrow\,
A\times_M A$ defined by $(a,v)\longmapsto (a,a+v)$ is an isomorphism.
In particular, for $x\in M$ the fiber $A_x$ is an
affine space over the vector space $V_x$.

Fix a theta characteristic $\Ohalf$ on $X$ -- in other words a
holomorphic line bundle $\Ohalf$ equipped with an isomorphism
$(\Ohalf)^{\ot 2}\cong \Omega_X$. If $X$ is compact of genus $g$,
there are $2^{2g}$ possible (distinct) choices. For
any $m\in {\mathbb Z}$, we will denote
$(\Ohalf)^{\ot m}$ by $\Omega_X^{\frac{m}{2}}$. The constructions
below are independent of the choice of $\Ohalf$ (see
\remref{indep of theta}).

\subsection{Kernel functions.}\label{derham-kernel}
Let $V$, $W$ be vector bundles over the Riemann surface
$X$, and $\on{Diff}^n(V,W)$ the sheaf of differential operators over
$X$ of order $n$ from $V$ to $W$.  Differential operators are a
bimodule over the sheaf $\Oo$ of holomorphic functions, via pre-- and
post--multiplication. We identify $\Oo_X$--bimodules with sheaves on
$X\times X$ (via the two pullback maps from $\Oo_X$ to
$\Oo_{\XX}$). This way (following Grothendieck and Sato) differential
operators are identified with ``integral kernels'' on $X\times X$ (see
e.g. \cite{BS,book}): there is a canonical isomorphism of
$\Oo_X$--bimodules (supported on the divisor $(n+1)\Delta$)
$$ \on{Diff}^n(V,W) = \frac{W\boxt V^\vee((n+1)\Delta)}{W\boxt
V^{\vee}} \, .$$ This is a coordinate--free reformulation of the
Cauchy integral formula: differential operators of order $n$ on
functions may be realized as kernel functions of the form
$\psi(z_1,z_2)dz_2$ with pole of order $(n+1)$ at $z_1=z_2$, via the
assignment
$$
f(z)\longmapsto \on{Res}_{z_1=z_2} f(z_1)\psi(z_1,z_2)dz_2\, .
$$ The
resulting differential operator depends only on the polar part of
$\psi$, equivalently on the restriction of $\psi$ in
$\Oo_X\boxt\Omega_X((n+1)\Delta)|_{(n+1)\Delta}$. For example, the de
Rham differential $d\,:\, {\mc O}_X\, \longrightarrow\, \Omega_X$ is
given by a holomorphic section $\mu_d$ of $\Omega_X\boxtimes \Omega_X
(2\Delta)$ over $2\Delta$. In local coordinates this
section is given
by $\dfrac{dz_1\boxt dz_2}{(z_1-z_2)^{2}}$, where
$z_i$ the pullback of a coordinate function $z$ on $X$
using the projection to the $i$th factor. 

\subsubsection{Connections as kernels.}\label{flat connections}

Let $\nu$ be an integer, and consider the line bundle
$${\mc M}_{\nu} = \Omega_X^{\frac{\nu}{2}}\boxt
\Omega_X^{\frac{\nu}{2}}(\nu\Delta)$$
over $X\times X$. As we have observed in \secref{derham-kernel},
the de Rham differential $d$ defines a section $\mu_d$ of
${\mc M}_2$ over $2\Delta$. In fact, for any $\nu
\in {\mathbb Z}$ there is a unique
trivialization $\mu_\nu$ of $\M_\nu|_{2\Delta}$ such that
\begin{enumerate}
\item{} $\M_\nu|_{\Delta}\cong \Oo_X$ is the natural trivialization
(defined by adjunction) (in other words, $\mu_\nu|_{\Delta}=1$);

\item{} the trivialization is symmetric, i.e., respects the
identification $\M_\nu\cong \sigma^*\M_\nu$ in the sense that
$\sigma^*\mu_\nu=(-1)^\nu \mu_\nu$.
\end{enumerate}
(Recall that $\sigma:X\times X\to X\times X$
is the interchange of factors.)
In particular note that $\mu_d\,=\, \mu_2$ and
$\mu_\nu\,=\, (\mu_1)^{\otimes\nu}$.
For a vector bundle $E$, denote by $\M_\nu(E)$ the vector bundle 
$$
\M_{\nu}(E)\,=\, E\boxt E^* \ot \M_\nu\, =\, 
(E\ot\Omega_X^{\frac{\nu}{2}})
\boxt(E\ot \Omega_X^{\frac{\nu}{2}})(\nu\Delta)
$$
on $X\times X$.

Consider the space $\Conn(E)$ of holomorphic connections on $E$. A
connection is given, following Grothendieck, by an isomorphism between
the two pullbacks $p_1^*E=E\boxt\Oo$ and $p_2^*E=\Oo\boxt E$
over $2\Delta$ (the
first--order infinitesimal neighborhood of the diagonal),
which restricts to the identity automorphism of $E$ on the diagonal.
In other words, a connection is determined by a section of
$\M_0(E)=E\boxt E^*$ on $2\Delta$ with ``symbol'' the identity
map $\text{Id}_E$ on the diagonal.

A connection on $E$ may also be described as a first--order
differential operator
$$
\nabla\,:\, E\,\longrightarrow\, E\otimes \Omega_X
$$
whose symbol is
the identity map $\text{Id}_E$. Thus $\nabla$ gives rise to a section of
$$\M_2(E)=(E\ot\Omega_X)\boxt (E^*\ot\Omega_X)(2\Delta)$$
on $2\Delta$
with biresidue the identity. These two formulations are related by
tensoring with the de Rham kernel $\mu_2=\mu_d$ trivializing $\M_2$ on
$2\Delta$. Similarly we can identify connections with sections of
$\M_\nu(E)|_{2\Delta}$ for any $\nu$. Note also that the
difference between any two connection kernels is a section of
$\Omega_X\ot\on{End}E$. Thus $\Conn(E)$ is an affine space for the
space $H^0(X,\, \Omega_X\ot\on{End}E)$ of endomorphism--valued
one--forms, or {\it Higgs fields}, on $E$.

Any holomorphic connection on a Riemann surface is flat (since there
are no nonzero holomorphic
two--forms on $X$.) This means that the identification between
nearby fibers of $E$ can be uniquely extended to an isomorphism
$p_1^*E\to p_2^*E$ to any order along $\Delta$ (in fact to local
trivializations in the analytic topology). Equivalently there is a
canonical extension from sections of $E\boxt E^*|_{2\Delta}$ which are
identity on the diagonal to sections $\kappa_\nu$ on $\nu\Delta$ for
any $\nu>0$, which in terms of a local flat basis of sections
$\{e_i\}$ with dual basis $\{e_i^*\}$, is given by
$$
\kappa_n\, =\, \sum e_i\boxt e_i^*\,\in\,
\Gamma((E\boxt E^*)|_{n\Delta})\, .
$$

In particular we obtain an isomorphism
$E\boxt E^*|_{n\Delta}\cong p_1^*\on{End}E|_{n\Delta}$.
In the language of kernels this map may be described as
the composition
$$E\boxt E^* \, \stackrel{\ot
\kappa_n^t}{\longrightarrow}\, \on{End}E\boxt \on{End}E\,
\stackrel{\on{tr}_E\ot\on{Id}}{\longrightarrow}\, \on\Oo
\boxt \on{End}E\, .
$$
Here $\kappa_n^t = \sigma^*\kappa_n \in \Gamma(E^*\boxt E|_{n\Delta})$
is the transpose of $\kappa_n$,
and $\on{tr}_E$ is the trace divided by the rank of $E$.

Note that this extension is
nonlinear with respect to the affine structure on $\Conn(E)$: it
involves solving the differential equation defining flat sections.

\subsection{Opers and kernel functions}\label{projective structures}

We would like to consider monic $n$th order differential operators
$$
L\, =\, \partial_t^n-q_1\partial_t^{n-1}-q_2\partial_t^{n-2}-
\cdots-q_{n-1}\partial_t -q_n
$$
on a Riemann surface $X$. To make this notion coordinate--independent,
we take $L:\Ac\to \Ac'$ to be a $n$th order operator between two
holomorphic line bundles, whose principal symbol is an isomorphism.
Since the symbol is a section of $\on{Hom}(\Ac,\Ac'\ot\Omega_X^{-n})$,
we must have $\Ac'\cong \Ac\ot\Omega_X^{\ot n}$.

It is
convenient to label the differential operator $L$ not by the line
bundle $\Ac$ but by its twist $\Ll=\Ac\ot\Omega_X^{\frac{n-1}{2}}$:

\subsubsection{Definition.} A $\text{GL}_n$--oper on $X$
consists of the data of a
line bundle $\Ll$ and a monic $n$th order differential operator
$$L\,\in\,
\Gamma(\on{Diff}^n(\Ac,\Ac\ot\Omega_X^n))\,=\,
\Gamma(\on{Diff}^n(\Ll\ot
\Omega_X^{\frac{1-n}{2}},\Ll\ot\Omega_X^{\frac{1+n}{2}}))
$$
over $X$ where
$\Ac=\Ll\ot\Omega_X^{\frac{1-n}{2}}$. The space of all
$\text{GL}_n$--opers
on $X$ is denoted by $\Op_n$, and opers for given $\Ll$ by
$\Op_n(\Ll)$.

\subsubsection{} It follows from the differential operator--kernel 
dictionary that
$\text{GL}_n$--opers for given $\Ll$ correspond to kernel
functions in $\M_{n+1}(\Ll)$ on $(n+1)\Delta$,
whose restriction to the
diagonal is the constant $1$ (by the trivialization
defined using adjunction).

Moreover, note that restricting the kernel function to $2\Delta$ we
obtain a section of $\M_{n+1}(\Ll)|_{2\Delta}$, which by \secref{flat
connections} defines a connection on $\Ll$. (This is the reason for
labeling opers by $\Ll$ rather than by $\Ac$.) Thus we have a
morphism $\Op_n(\Ll)\longrightarrow \Conn(\Ll)$. In particular, for
$\Ll=\Oo$, we can look for opers which induce the trivial connection
on $\Oo$, so that the associated kernel on $2\Delta$ is the de Rham
kernel $\mu_{n+1}$. In terms of differential operators, the induced
connection (restriction to $2\Delta$) is determined by the {\em
subprincipal symbol} $q_1$. Thus we are considering differential
operators $L$ of the form
$$L\,=\,\partial_t^n-q_2\partial_t^{n-2}-\cdots-q_n.
$$
(Conversely the vanishing of the subprincipal symbol forces
$\Ll$ and the associated connection to be trivial.)

\subsubsection{Definition.} A
$\text{SL}_n$--oper on $X$ (for fixed theta characteristic $\Ohalf$)
is a monic $n$th order differential operator
$$
L\,\in \,\Gamma(\on{Diff}^{n}
(\Omega_X^{\frac{1-n}{2}},\Omega_X^{\frac{1+n}{2}}))
$$ 
with vanishing subprincipal symbol. Equivalently, $L$ is defined by a
section of $\M_{n+1}$ on $(n+1)\Delta$, whose restriction to $2\Delta$
agrees with $\mu_{n+1}$. The space of $\text{SL}_n$--opers (for fixed
$\Ohalf$) is denoted by $\Op_n^{\circ}$.

\subsubsection{Remark.}\label{indep of theta}
The restriction of the bundles $\M_n$ to any neighborhood $k\Delta$
are independent of the choice of theta characteristic $\Ohalf.$ This
follows from the fact that the ratio of two theta characteristic is a
bundle of order two, $\Ll^2=\Oo_X$, and so carries a canonical flat
connection (inducing the trivial connection on ${\Oo}_X$), which gives
rise to a trivialization of $\Ll\boxt \Ll^*$ on $n\Delta$ for any $n$.
This may also be seen from the universal form of the transition
functions defining $\Omega_X^{\frac{\nu}{2}}\boxt
\Omega_X^{\frac{\nu}{2}}|_{(n+1)\Delta}$ -- in fact these transition
functions make sense for an arbitrary {\it complex} number $\nu $,
since the Taylor expansion of an expression $\dfrac{dz_1^\nu\boxt
dz_2^\nu}{(z_1-z_2)^{2\nu}}$ in terms of a new coordinate $w=w(z)$ has
coefficients which are polynomials (with integer coefficients) in $\nu
$. In other words, all of these bundles are attached to natural
representations of the group of formal changes of coordinates on $X$
(see \cite[7.2]{book}).

Thus the spaces of opers for different choices of $\Ohalf$ are all
isomorphic. Alternatively, one can define $\text{PSL}_n$--opers and then
identify $\text{SL}_n$--opers with pairs consisting of a
$\text{PSL}_n$--oper and a theta characteristic (\cite{opers}).

\subsubsection{Example.} On $\Pp^1$, there is a unique
$\text{SL}_n$--oper for every $n$ (here
$\Omega_{\Pp^1}^{\frac{1}{2}}=\Oo_{\Pp^1}(-1)$ is the unique theta
characteristic). It is defined by the kernel function
\begin{equation}\label{BR formula}
\gamma_\nu\,=\, \frac{dz^{\frac{n+1}{2}}\boxt
dz^{\frac{n+1}{2}}}{(z_1-z_2)^{n+1}}
\end{equation} 
on $(n+1)$, where $z$ is the natural coordinate function on
$\C \, \subset\, \C\cup\{\infty\} = \Pp^1$. This $\gamma_\nu$
is a holomorphic section over $\Pp^1\times \Pp^1$, and it
invariant under the diagonal action of $\text{PSL}_2$ on
$\Pp^1\times \Pp^1$.

\subsubsection{Lemma.}\label{GL_n to SL_n} 
There is a canonical isomorphism
$\Op_n(\Ll)\, =\, \Conn(\Ll)\times\Op_n^\circ$.

\subsubsection{Proof.} An oper $L\in \Op_n(\Ll)$ defines a connection
on $\Ll$ as above. Solving the connection defines a trivialization
$\kappa_{n+1}$ of $\Ll\boxt\Ll^*$ on $(n+1)\Delta$. This
trivialization gives an
isomorphism $\M_{n+1}(\Ll)\cong \M_{n+1}$, which sends $L$ to an
$\text{SL}_n$--oper $L'$. The kernel of $L'$ is explicitly given by
$\kappa_{n+1}\inv$ times the kernel of $L$, from which it is obvious
that the restriction to $2\Delta$ is indeed $\mu_{n+1}$.

\subsubsection{Projective Structures.}

A {\it projective structure} on a Riemann surface $X$ (see
\cite{Gu}, \cite{De}) is an equivalence class of atlases $\{U_{\alpha},
{\vphi}_{\alpha}\}_{\alpha\in I}$ on $X$, where ${\vphi}_{\alpha}$ is
a holomorphic embedding of the open set $U_{\alpha}$ in $\Pp^1$, so
that the transition maps ${\vphi}_{\beta}\circ {\vphi}^{-1}_{\alpha}$
are M\"obius (or fractional linear) transformations (elements of
$\text{PSL}_2\C$). The space of projective structures will be denoted
$\Proj$. A projective structure on $X$ allows us to pull back any
$\text{PSL}_2$--invariant construction from $\Pp^1$ to $X$. In
particular, we may pull back the $\text{SL}_n$ opers on $\Pp^1$ (and
their kernel functions $\gamma_{n+1}$) to define $\text{SL}_n$--opers
for every $n$, or equivalently monic differential operators $D_{n}$
with vanishing subprincipal symbol. The symbol of each such operator
is the constant function $1$. The operator $D_0$ is the identity
automorphism of $\Ohalf$. The operator $D_1$ is the exterior
derivative $d\, :\, {\mc O}_X \, \longrightarrow\, \Omega_X$. The
operator $D_2 \,\in\, \Gamma(\on{Diff}^2(\Omega_X^{-\frac{1}{2}},
\Omega_X^{\frac{3}{2}}))$ over $X$ is the {\em Sturm--Liouville
operator}, or {\em projective connection}, associated with a
projective structure. Thus $D_2$ is the differential operator which
in {\em projective} local coordinates has the form $\pa_t^2$.

The projective structure can be recovered from the associated
projective connection, setting up a bijection $\Proj\cong
\Op_2^\circ$: the projective atlases are defined by the ratios of any
two local linearly independent solutions of the Sturm--Liouville
operator $D_2$.

\subsection{Opers as connections.}\label{opers as connections}
Opers have an interpretation in terms of vector bundles with
connection, which also enables the generalization from $\text{GL}_n$
to an arbitrary reductive group. This observation and its current
formulation are due to Drinfeld--Sokolov \cite{DS} and
Beilinson--Drinfeld \cite{opers}, respectively. Recall that the study
of the differential operator
$$
L\,=\,\partial_t^n-q_1\partial_t^{n-1}-q_2\partial_t^{n-2}-\cdots-q_n
$$
is equivalent to that of the system of $n$ first--order equations
which can be written in terms of the first--order matrix operator
$$\partial_t-\left( \begin{array}{ccccc}
q_1&q_2&q_3&\cdots&q_n\\
1&0&0&\cdots&0\\
0&1&0&\cdots&0\\
\vdots&\ddots&\ddots&\cdots&\vdots\\
0&0&\cdots&1&0
\end{array}\right).$$
Now suppose $L$ is a
$\text{GL}_n$--oper on $X$ for the line bundle $\Ac$.
It is not hard to see that the above first--order systems patch
together to define a connection $\nabla:F\to F\ot \Omega_X$ on a rank
$n$ vector bundle $F$, which carries a filtration $0=F_0\subset
F_1\subset \cdots\subset F_n=F$, with $F_1\cong {\mc A}$. The key
features of the above matrix system are the appearance of zeros
beneath the subdiagonal (Griffiths transversality) and $1$s on the
subdiagonal (nondegeneracy). Locally, the bundle and flag
$(F,F_{\bullet})$ admit a unique trivialization so that the connection
has the above form. Moreover for the connection to be an
$\text{SL}_n$--connection (so that the determinant line bundle and its
connection are trivial) the subprincipal symbol $q_1$ must vanish, so
that we obtain $\text{SL}_n$--opers.

\subsubsection{Proposition.}\label{BD prop}
(\cite{opers}) $\text{GL}_n$--opers on $X$ correspond 
canonically to the data of a
rank $n$ vector bundle $F$, equipped with a flag
$$0\subset F_1\subset\cdots F_{n-1}\subset F_n=F,$$ and a connection
$\nabla$, satisfying
\begin{enumerate}
\item[$\bullet$] $\nabla(F_i) \subset F_{i+1} \otimes \Omega_X$.
\item[$\bullet$] The induced maps $F_i/F_{i-1}\to (F_{i+1}/F_i)
\otimes \Omega_X$ are isomorphisms for all $i$.
\end{enumerate}
$\text{SL}_n$--opers are $\text{GL}_n$--opers for which the determinant
line bundle of the flat vector bundle $(F,\nabla)$ is trivial.

\subsubsection{} In fact, the transversality condition on the
connection is sufficiently
rigid to force the underlying vector bundle $F$ to be the
$(n-1)$st jet
bundle $F\cong J^{n-1}(F/F_{n-1})$, with its canonical
filtration. 

Note that from the connection point of view, the extension of 
projective connections to $n$th order operators 
is simply the operation of inducing an
$\text{SL}_n$--oper from a $\text{SL}_2$--oper 
by taking the associated bundle for the
$(n-1)$st symmetric power representation of
$\text{SL}_2$ into $\text{SL}_n$.

\subsection{The Hitchin base.}
An important non-obvious feature of opers (for fixed $\Ll$) is that
they form an affine space over the Hitchin base space of $X$
(\cite{Hitchin}),
$$
\on{Hitch}_n(X)\, =\, \Gamma(\Omega_X)\oplus
\Gamma(\Omega_X^{2})\oplus\cdots \oplus
\Gamma(\Omega_X^n)\, .
$$
This generalizes the statement that projective
structures form an affine space over quadratic differentials
$\Gamma(\Omega_X^2)$.

\subsubsection{Proposition.}\label{Hitchin opers}
(\cite{opers}) There is a canonical isomorphism
$$\Conn(\Ll)\times \Proj\times
\on{Hitch}_n^{>2}\,\longrightarrow\, \Op_n(\Ll)\, .$$

\subsubsection{Remarks on proof.} 
Geometrically, the proposition is an expression of the fact that the
tangent bundle to $\Pp^{n-1}$ restricted to the rational normal curve
splits canonically (i.e., $\text{PSL}_2$--equivariantly) into a sum of
line bundles. Namely, the stabilizer in $\text{PSL}_2$ of a point on
the rational normal curve (which is isomorphic to upper triangular
matrices $B_0$) acts on the tangent space at that point through its
$\Cx$--quotient. In fact since a $\text{SL}_n$--oper naturally gives
rise to a projective structure (on restriction to $3\Delta$), the
proposition reduces to this fact since any infinitesimal (or
complex-local) statement on $\Pp^1$ which is
$\text{PSL}_2$--equivariant generalizes to any Riemann surface with
projective structure. In particular, given an oper $F$ induced from a
$\text{SL}_2$--oper, one identifies a subbundle $\V\cong\bigoplus_1^n
\Gamma(\Omega_X^i)$ inside the Higgs fields
$\Gamma(\on{End}F\ot\Omega_X)$.
Addition of sections from $\V$ acts transitively on oper
connections on $F$ -- in particular the action of $\Gamma(\Omega_X)$
changes the connection, while that of $\Gamma(\Omega_X^2)$ changes
the projective structure.

\subsubsection{Remark.}\label{dim of opers} 
It follows from the proposition that the dimensions of
$\Op_n$ and $\Op_n^\circ$
on a compact Riemann surface $X$ of genus $g$ are $(g-1)(n^2-1)+g$ and
$(g-1)(n^2-1)$ respectively.

\subsection{Shifted opers.}\label{shifted opers}
The projection from $\Op_n^\circ$ to $\Proj$ may be described 
conveniently using kernels.
Given an $\text{SL}_n$--oper with kernel $s\in
\Gamma(\M_{{n+1}}|_{(n+1)\Delta})$, 
its restriction to $3\Delta$ defines an element of the space
$$
\Proj(k)\,=\,\{s\in
\Gamma(\M_{n+1}|_{3\Delta}) \; \vert\; s|_{2\Delta}=\mu_{n+1}\}\, .
$$
These ``shifted'' projective kernels are however naturally identified
with projective structures.
Note that the difference between any two sections of $\Proj(k)$ vanishes
on $2\Delta$, and so may be identified with a section of
$\M_k(-2\Delta)|_{\Delta}\cong \Omega_X^{\ot 2}$, that is, a quadratic
differential on $X$. It follows that
$\Proj(k)$ is a torsor for the quadratic differentials
$\Gamma(\Omega_X^{\ot 2})$ on $X$.
Recall that we may rescale the torsor structure on a fixed affine
bundle by any scalar $\la\in\C^{\times}$, keeping the manifold
$\pi:A\to M$ the same but making $v\in V$ act by $\la\cdot v$.

\subsubsection{Lemma.}\label{rescaling torsors} The spaces
$\Proj(k)$ for $k\neq 0$ are all isomorphic
(with rescaled torsor structure over quadratic differentials).

\subsubsection{Proof.} The $k$--th power map
$\rho\mapsto\rho^{\ot n}|_{3\Delta}$ 
identifies the sheaves $\Proj(1)$ and $\Proj(k)$. On affine
structures, this has the effect of rescaling by $k$:
$$
(\rho+q)^{\ot k}|_{3\Delta}\,=\,(\rho^{\ot k}+k\cdot q)|_{3\Delta}
$$
(since all higher terms vanish on $3\Delta$). Note that
multiplication by $k$ sends $\Proj(1)$ isomorphically to sections of
$\M_1|_{3\Delta}$ whose restriction to $2\Delta$ is $k\mu_1$, and has
the same effect on torsor structures.

\subsubsection{}
It follows that the restriction of a
$\text{SL}_n$--oper to $3\Delta$ may be
naturally identified with a projective structure on $X$. (For
$\text{GL}_n$--opers, we must twist by the connection on $\Ll$
given by the restriction to $2\Delta.$)

One naturally encounters other kernel realizations of the spaces of
opers:

\subsubsection{Definition.} The space 
of {\em $k$--shifted opers} on $\Ll$ is defined to be the space 
$$
\Op_n(\Ll)(k)\,=\,\{s\in\Gamma(\M_k(\Ll)|_{(n+1)\Delta}) \; \vert
\; s|_{\Delta}=1\}\, .
$$

\subsubsection{}
(Note that for $k>n$ shifted opers form a quotient of
$\Op_{k-1}^{\circ}$ by differential operators of order $k-n-2$.)
The restriction of a shifted oper to $3\Delta$ defines an element of
$\Proj(k)$, hence a projective
structure on $X$. It follows that we may identify all shifted
opers canonically with honest opers. Explicitly, this is done by
tensoring with the kernels $\gamma_{{n+1}-k}$ obtained from the
projective structure by \eqref{BR formula}, which trivializes
$\M_{{n+1}-k}$ to any order near the diagonal.

\subsection{Projective kernels and projective connections.}\label{Bergman} 
Many projective connections arising in Riemann surface theory arise
naturally from {\em projective kernels}, or global bidifferentials
of the second kind, with biresidue one:

\subsubsection{Definition.} 
\begin{enumerate}
\item A $\text{SL}_n$--oper kernel on $X$
is a global section $s\in \Ho(\XX, \M_n)$ with $s|_{2\Delta}=\mu_n$.
The space of oper kernels is denoted by $\Kern_n$.

\item A projective kernel on $X$ is a symmetric $\text{SL}_2$--oper
kernel. In other words, a bidifferential
$$
\omega\,\in\, \Ho(\XX,\Omega_X\boxt\Omega_X(2\Delta))^{\text{sym}}
$$
with biresidue
one on $\XX$. The space of projective kernels is denoted by
$\Kern_2^{sym}\subset\Kern_2$. 
\end{enumerate} 
By $\Ho(\XX,\Omega_X\boxt\Omega_X(2\Delta))^{\text{sym}}$ we mean
section invariant under the involution $\sigma$.

\subsubsection{} The
difference between two projective kernels is a holomorphic symmetric
bidifferential on $X$, so that $\Kern_2^{sym}$ is an affine space for
$\on{Sym}^2\Ho(X,\Omega)$, of dimension $\binom{g}{2}$. Restriction to
$3\Delta$ defines a map $\Kern_2^{sym}\to\Proj_2$, which is surjective
for $X$ non--hyperelliptic.

A key r\^ole of the spaces
$\Proj$ and $\Kern_2^{sym}$ is in relation to moduli
spaces. Namely, $\Ho(X,\Omega^{\ot 2}_X)$ is the cotangent space to the
moduli space of curves (or Teichm\"uller space) at $X$, and
$\on{Sym}^2\Ho(X,\Omega_X)$ is the cotangent space to the moduli of
abelian varieties (or Siegel upper half space) at the Jacobian
$\Jac$ of $X$. The spaces $\Proj$ and $\Kern_2^{sym}$ are naturally
identified as the fibers, at $X$ and $\Jac$ respectively, of the
space of connections on a natural (Hodge or theta) line bundles on the
respective moduli spaces (in other words the fibers of appropriate 
twisted cotangent bundles). (See \cite{Tyurin, part 1}.) 

An important example of a projective kernel is the Bergman kernel
$\omega_B$. Let $\omega_i$ ($i=1,\dots,g$) be the normalized basis of
holomorphic differentials on $X$, with respect to a normalized homology
basis $A_i,B_j$,
and
$\frac{\partial}{\partial z_i}$ the dual basis of vector fields
on the Jacobian. 
The Bergman kernel is 
characterized by having vanishing $A$ periods, and the forms
$\omega_i$ as its $B$--periods.

\subsection{Matrix opers.}\label{mopers}
In this section we describe a matrix version of opers. Thus we consider
$n$th order differential
operators with matrix coefficients,
$$
L\,=\,\partial_t^n-q_1\partial_t^{n-1}-q_2\partial_t^{n-2}-\cdots-q_n
$$
where the $q_{n+1}$ are now $k$ by $k$ matrices, acting on $\C^r$. 
Let $L:E_1\to E_2$ be a $n$th
order differential operator between vector bundles $E_1,E_2$ on $X$,
whose symbol is an isomorphism, so that $E_2\cong E_1\ot \Omega_X^{\ot
n}$.

\subsubsection{Definition.} A $n$th order matrix oper on $E$ is an
$n$th order differential operator
$L\,\in\,\Gamma(\on{Diff}^n(E\ot\Omega_X^{\frac{1-n}{2}},
E\ot\Omega_X^{\frac{1+n}{2}}))$ over $X$
with principal symbol the identity
$\on{Id}_E$. The space of matrix opers on $E$ is denoted by $\Mop_n[E]$.

\subsubsection{}
We may follow the same procedure as in \secref{opers as connections}
to describe $n$th order matrix opers $L$ by first order matrix
systems, now of rank $nk$. The resulting connections were called {\it
coupled connections} in \cite{coupled connections}. We have the
following statement (see
\cite{coupled connections} for more details):

\subsubsection{Proposition.}\label{generalized opers} 
Let $E$ be a vector bundle. There is a natural identification between
matrix opers $L:E\to E\ot \Omega_X^{\ot n}$ of order $n$, and vector
bundles $F$ equipped with a filtration $0=F_0\subset F_1\subset
\cdots\subset F_n=F$, with $F_1\cong E\ot \Omega_X^{\frac{1-n}{2}}$,
and a connection $\nabla:F\to F\ot \Omega_X$ satisfying the two
conditions
\begin{enumerate}
\item[$\bullet$] $\nabla:F_\nu\to F_{\nu+1}\ot \Omega_X$ (Griffiths
transversality), and
\item[$\bullet$] the homomorphism $F_\nu/F_{\nu-1}\to
F_{\nu+1}/F_\nu\ot \Omega_X$ induced by $\nabla$ is an isomorphism
for all $\nu$.
\end{enumerate}

\subsubsection{Proof.}
Recall that a $n$th order operator
$L\in\Gamma(\on{Diff}^n(E_1,E_1\ot\Omega_X^n))$
over $X$ is a homomorphism from $J^n(E_1)$ to
$E_1\ot\Omega_X^n$. This is equivalent to a splitting of the jet
sequence
$$0\to E_1\ot\Omega_X^n \to J^n(E_1)\to J^{n-1} (E_1)\to 0,$$
and thus
to a lift $J^{n-1}(E_1)$ to $J^n(E_1)$. However there is a natural
homomorphism $J^n(E_1)\to J^1(J^{n-1}(E_1))$ for any bundle. Thus we
have constructed a lifting from $J^{n-1}(E_1)$ to its sheaf of
$1$--jets, in other words a connection on $J^{n-1}(E_1)$. The strict
Griffiths transversality with respect to the natural connection on
$J^{n-1}(E)$ follows automatically.

In the reverse direction, given a filtered vector bundle $F$ and a
connection $\nabla$ on $F$ as above,
consider the homomorphism
\begin{equation*}\label{split}
\psi_k\, :\, F\, \longrightarrow \,
J^k(F)\,\longrightarrow\,J^k(F/F_{n-1})
\end{equation*} 
 where the first arrow is the flat extension map given by the
 connection and the second is the projection. The transversality
 condition ensures that $\psi_{n-1}$ is an isomorphism. Therefore,
 $\psi_{n}\circ\psi_{n-1}\inv:J^{n-1}(F/F_{n-1})\to J^n(F/F_{n-1})$
 gives a splitting of the jet sequence as above, in other words, a
 differential operator as desired.

\subsubsection{Developing maps.}
A geometric description of coupled connections $\nabla$ as in
\propref{generalized opers}, generalizing the description of projective
structures via period maps, is given in \cite{coupled connections}.
Namely consider the Grassmannian bundle $G_k(F)$ of $k$--dimensional
subspaces of $F$. This inherits a connection from $\nabla$ and a
section from $F_1$, which is nowhere flat. It follows that on simply
connected opens (or on the universal cover of $X$) we obtain period
maps to $G_k(\C^{nk})$ using the connection to trivialize $G_k(F)$ and
the section to map. These period maps satisfy natural nondegeneracy
conditions. Conversely such nondegenerate period maps with transitions
coming from the action of $\text{GL}_{nk}$ on $G_k$ give rise to coupled
connections. 

\subsubsection{Decomposition of matrix opers.}\label{decomposition}
Matrix opers of order $n$ on $E$ correspond to kernel functions
in $\M_{{n+1}}(E)|_{(n+1)\Delta}$, that is a section of
$\M_{{n+1}}(E)$ over $(n+1)\Delta$,
whose restriction to the diagonal is
$$
\on{Id}_E\,\in\,\on{End}E\,\cong\, \M_{n+1}(E)|_{\Delta}\, .
$$
For example, if
$E=\Ll$ is a line bundle, then matrix opers and $\text{GL}_n$--opers for
$\Ll$ are the same. It follows that by restriction to $2\Delta$, a
matrix oper defines a section of $\M_{{n+1}}(E)|_{2\Delta}$ with residue
$\on{Id}_E$ and thus a flat connection on $E$ (\secref{flat
connections}). So there is a canonical projection
$\Mop_n(E)\to\Conn(E)$.

The induced connection on $E$ allows us to identify
$\M_{n+1}(E)$ with $\M_{{n+1}}\ot p_1^*\on{End} E$ to any order
near $\Delta$. Thus, if $p\in \C[{\mf gl}_n]^{\text{GL}_n}$ is an
invariant polynomial on matrices (i.e., a coefficient of the
characteristic polynomial) we obtain a map 
$$
p_*\,: \,\Mop_n(E)\,\longrightarrow\, \Op_n^{\circ}
$$ by applying $p$ to
$\on{End}E$ and identifying the resulting shifted oper with an oper.

Together with
\propref{Hitchin opers}, this gives a very simple description of
matrix opers. Let
$$\on{Hitch}_n^{>1}(E)^\circ\,=\,\bigoplus_{i=2}^n
\Gamma(\on{End}^\circ E\ot\Omega_X^{\ot i}),$$
the space of {\em traceless} $\on{End}E$--valued polydifferentials.

\subsubsection{Proposition.}\label{Hitchin mopers}
There is a canonical isomorphism
$$
\Conn(E)\times\Op_n^\circ\times 
\on{Hitch}_n^{>1}(E)^{\circ}\,\longrightarrow\, \Mop_n(E)\, .
$$

\subsubsection{Proof.} 
We describe the isomorphism in the languages of kernels and of coupled
connections.

We define the map $\Conn(E)\times\Op_n^\circ\to \Mop_n(E)$ using the
tensor decomposition of $\M_n(E)$ by taking the tensor product of
sections. It follows that the decomposition (\propref{Hitchin opers})
of sections of $\M_{n+1}$ gives rise to a direct sum decomposition of
sections of this tensor product. We identify $\Op_n^\circ$ with the
scalar endomorphisms, thereby obtaining the proposition. The
projection back to $\text{SL}_n$--opers is given by the induced map
$\on{tr}_*$ for $p(A)=\on{tr}(A)/\on{rk}(E)$ above.

Viewing an oper through the corresponding flat bundle
$(F,F_{\bullet},\nabla)$, where $F_{\bullet}$ is a filtration of
subbundles of $F$,
we may take the tensor product of vector bundles
$E\ot F$, with its induced filtration and connection. The result is a
coupled connection, which we consider as a matrix oper on $E$. Inside
the space of Higgs fields on $E\ot F$ compatible with the filtration
we find the tensor product of $\on{End}E$ with the space $\V$ of Higgs
fields from \propref{Hitchin opers}, so we can modify the coupled
connection by $\on{End}E$--valued polydifferentials. Again one checks
this gives a bijective parametrization of coupled connections.

To see the compatibility of the constructions, note that a vector
bundle with connection is canonically trivialized to any order near a
point $x\in X$, up to a constant matrix (the change of trivialization
of its fiber at $x$). Hence the compatibility reduces to the
(equivariant) compatibility in the case of the trivial bundle, which
is obvious.

\subsubsection{The determinant.} 
The determinant map for matrix opers may also be described directly,
without solving the associated connection.
Let $s$ be a section of $(E_1\boxt E_2)\ot\Ll$, where $E_1, E_2$ are
vector bundles on $X$ of the same rank $k$ and $\Ll$ is a line bundle
on $X\times X$. Then we may define the determinant section $\det
s=\wedge^k s$ of $(\det E_1\boxt\det E_2) \ot \Ll^k$ (e.g. consider
$s$ as a homomorphism from $p_1^*E_1$ to $p_2^*E_2 \ot\Ll$ of rank $k$
vector bundles and take its determinant).

If $s\in E\boxt E^*|_{2\Delta}$ is a connection on $E$, then $\det
s\,\in\, \Gamma(\det E\boxt \det E^*|_{2\Delta})$ is the determinant
connection on $\det E$. More generally, 
the determinant defines a canonical map
$$
\det\, :\, \Mop_n(E)\,\longrightarrow\, \Op_n(\det E)
$$
Namely, the determinant of $s\,
\in\,\Gamma( \M_n(E)|_{k\Delta})$ defines a section
$$
\det s\,\in\, \Gamma(\M_{n\on{rk} E}(\det E)|_{k\Delta})
$$
which is the identity on the diagonal, i.e., a shifted oper, and
which we identify with an (unshifted) oper as in \secref{shifted
  opers}.  There is a commutative diagram
\begin{equation*}
\begin{array}{ccc}
\Mop_n(E) & \longrightarrow & \Conn(E)\times
\Gamma(\M_n\ot p_1^*\on{End}E) \\
\Big\downarrow&&\Big\downarrow\\
\Op_n(\det E) & \longrightarrow & \Conn(\det E)\times\Op_n^{\circ}
\end{array}
\end{equation*}
where the horizontal arrows are given by trivializing $E,\det E$ using
the connection, and the vertical arrows are the determinant maps on
kernels and on endomorphisms. This identifies the determinant map for matrix
opers above with the determinant of the associated kernel.

\subsection{Extended Connections.}
The splitting in \propref{Hitchin mopers} picks out a particularly
interesting subspace $\Conn(E)\times\Proj$ of matrix opers on $E$, the
{\em extended connections} on $E$. In fact extended connections most
naturally appear as a {\em quotient} of matrix opers. Their r\^ole is
as an affine space for the cotangent space of the moduli of the pair
$(X,E)$. As such they do not split as a product: the splitting in
\propref{Hitchin mopers} is nonlinear (since it involves solving the
connection to some order), and in fact a deformation the quadratic part
of the Hitchin map.

\subsubsection{Definition.}\label{ex.con.def.}
The space $\excon_{n+1}(E)$ of extended connections on $E$ is the
space of monic sections of the quotient of $\M_{n+1}(E)|_{3\Delta}$ by
the subsheaf of sections vanishing on $2\Delta$ and with vanishing
trace on $3\Delta$.

\subsubsection{}
Here trace refers to the composition
$$
\M_{n+1}(E)(-2\Delta)\, \longrightarrow \,\Omega^{\otimes
2}_X\otimes\on{End}E\,\longrightarrow \,\Omega_X^2\, ,
$$
and monic sections
are sections restricting to $\on{Id}_E$ on the diagonal. Thus we have
modified $\M_{n+1}(E)|_{3\Delta}$ by forgetting all but the trace of
the lowest--order term.

It follows that restriction to $2\Delta$ makes $\excon_{n+1}(E)$ an
affine bundle for quadratic differentials $H^0(X,\Omega_X^2)$ over
$\Conn(E)$. Consider the space of {\em extended Higgs fields}
$$
\exhiggs(E)\,=\,\{s\in \Gamma(\M_{n+1}(E)|_{2\Delta}) \; \vert \;
s|_{\Delta}=0\}/\Gamma (\M_{n+1}(E)^{\circ}|_{2\Delta})\, .
$$
Note that this space $\exhiggs(E)$ is
independent of ${n+1}$ since $\M_{n+1}|_{2\Delta}$ is canonically
trivialized. The space of extended connections is clearly a torsor over
extended Higgs fields.
The importance of the latter is as the cotangent
space at $(X,E)$ to the moduli of pairs of Riemann surfaces and vector
bundles. They form an extension
\begin{equation*}
0\longrightarrow H^0(X,\Omega^{2}_X)\longrightarrow \exhiggs(E)
\longrightarrow  H^0(X,\on{End}E\ot \Omega_X)\longrightarrow  0
\end{equation*}
of Higgs fields on $E$ by quadratic differentials. It is proven in
\cite{part 1} that the torsors $\excon_{n+1}(E)$ over $\exhiggs(E)$
for varying $X,E$ form a twisted cotangent bundle over the moduli
space: in particular there is an isomorphism
$$
\excon_1\, \cong \,\Conn(\Theta)
$$
with the affine bundle of connections
on the theta line bundle over the moduli space.

\subsubsection{Lemma.} For every $n\in\Z$, 
the map
$$
\Conn(E)\times\Proj_{n+1}\,\longrightarrow \,
\Gamma(\M_{{n+1}}(E)|_{3\Delta})\,
\longrightarrow \, \Gamma(\M_{{n+1}}|_{3\Delta}/\M_{{n+1}}(E)^{\circ})
$$
defines an isomorphism
$\Conn\times\Proj\to\excon_{n+1}(E)$ and thereby lifts the latter to
$\Mop_2(E)({n+1})$ (and hence $\Mop_n$ for every $n$).

\subsubsection{The deformed quadratic Hitchin map.} 
The projection
$$\excon_{n+1}(E)\,\longrightarrow \, \Proj$$ of extended connections
back to projective structures may be described in several
ways. Following \secref{decomposition}, it is given by sending
$\M_{{n+1}}(E)|_{3\Delta}\to \M_{{n+1}}|_{3\Delta}\ot\on{End}E$ via
the connection and thence to $\Proj({n+1})$ via trace of sections.
Alternatively, it can be deduced from the determinant map
$$
\det\,:\,\M_{n+1}(E)|_{3\Delta}\,\longrightarrow \, \M_{k(n+1)}(\det
E)|_{3\Delta}\, .
$$
We identify the resulting $\text{GL}_2$--oper with an element of
$\Proj(k(n+1))$ by tensoring with the connection of $\det E$ as in
\lemref{GL_n to SL_n}, and thence with an element of $\Proj({n+1})$ by
\lemref{rescaling torsors}. (This agrees with the trace map since we
are restricting to $3\Delta$, thereby keeping only the leading term of
the determinant.)

Another description of the projection is given by taking the trace of
the square of the kernel. More precisely, for $s\in
\Gamma(\M_{n+1}(E)|_{3\Delta})$, its transpose $s^t=\sigma^*s\in
\M_{n+1}(E^*)|_{3\Delta}$, so that the tensor product lives in
$$
s\ot s^t\,\in\, (\on{End}E\boxt\on{End}E)\ot \M_{2(n+1)}
$$
over $3\Delta$. We apply trace to both factors, obtaining
$$
S(s)\,=\,\on{tr}_E\boxt\on{tr}_E (s\ot s^t)\,\in\,
\Gamma(\M_{2(n+1)}|_{3\Delta})
$$
which is monic if $s$ is. To compare
this with the other constructions, suppose $\rho
\in\M_{n+1}|_{3\Delta}$ is a projective structure, $\nabla$ is a
connection and $\kappa\in E\boxt E^*|_{3\Delta}$ is the corresponding
kernel function giving the isomorphism $p_2^*E\to p_1^*E$. Note that
$$
(\on{Id}\ot\on{tr}_E)(\kappa\ot\kappa^t) \,=\,\on{Id}_E\boxt
1\in\on{End}E\boxt \Oo\, ,
$$
simply expressing the fact that $\kappa^t$ is
the flat kernel for the inverse map $p_1^*E\to p_2^*E$. It follows
that $S(\rho\ot\kappa)=\rho$, so that $S$ is indeed the projection
back on projective structures.

This description of the determinant map for extended connections
presents it as a deformation of the quadratic Hitchin map. Namely let
$$
\excon_{n+1}^\la(E)\,=\,\{s\in \Gamma(\M_{n+1}(E)|_{3\Delta}) \; \vert
\; s|_{\Delta}=\lambda\on{Id}\}/\Gamma(\M_{n+1}(E)^{\circ}|_{3\Delta})
$$
be the family deforming extended connections to extended Higgs fields.

\subsubsection{Proposition.}
The determinant map $\excon_{n+1}(E)\to \Proj({n+1})$ deforms to a map
$$
\excon_{n+1}^\la(E)\,\longrightarrow \, \Proj(\lambda(n+1))
$$
(for $\la\in\C$),
which for $\lambda=0$ factors through the quadratic Hitchin map
$$
\exhiggs(E)\,\longrightarrow \, \Gamma(\Omega_X\ot\on{End}E) \,
\longrightarrow \,\Gamma(\Omega_X^2)\,=\,\Proj(0)\, ,
$$
sending $\eta\in \Gamma(\Omega_X\ot\on{End}E)$ to $\on{tr}_E(\eta^2)$.

\subsubsection{Proof.}
If $s|_{\Delta}=\lambda\on{Id}_E$ then
$S(s)|_{2\Delta}=\lambda^2\mu_{2(n+1)}$ (by symmetry with respect to
transposition of factors). For $\la\neq 0$ the space of such kernels
is isomorphic (by rescaling and taking square--root,
\lemref{rescaling torsors}) with projective structures. In fact the
resulting map $\excon_{n+1}^\lambda(E)\to\Proj(\lambda(n+1))$ is a
morphism of torsors for quadratic differentials (the square root
$\Proj(2\lambda(n+1))\to\Proj(\lambda(n+1))$ compensates for the
quadratic expression $s\ot s^t$.) This map clearly descends to
$\excon_{n+1}^\la(E)$. On the other hand, for $\lambda=0$, we obtain a
quadratic differential, realized as a section of $\M_{n+1}|_{3\Delta}$
vanishing on $2\Delta$. This quadratic differential depends only on
the Higgs field $\eta$ underlying the extended Higgs field $s$, and
equals $\on{tr}_E(\eta^2)$ (the first trace squares $\eta$ by
contracting indices, while the other trace takes trace of the
resulting matrix).

\section{The Klein and Wirtinger maps.}\label{Klein and Wirtinger}

Let $\Mf_X(n)$ denote the moduli space of semistable vector bundles
over $X$ of rank $n$ and Euler characteristic $0$. It is known that
$\Mf_X(n)$ is an irreducible normal projective variety of
dimension $(g-1)(n^2-1)+g$.
In particular, $\Mf_X(1)=\Pg$, the moduli of degree $g-1$ line bundles.
Let $\Mf_X(n)_0$ denote the moduli space of semistable vector bundles
of rank $n$ and degree $0$. The chosen theta characteristic
$\Ohalf$ gives an isomorphism 
$$\Mf_X(n)\longrightarrow \Mf_X(n)_0,\hskip.3in E\longmapsto
E_0=E\ot\Omega_X^{-\frac{1}{2}}$$
(since tensoring by a line bundle
preserves semistability).

The determinant map $E\mapsto \det E$ sends $\Mf_X(n)$ to
$\on{Pic}_X^{n(g-1)}$. We may identify a closed subvariety
$$
\N_X(n)\,=\,\on{det}^{-1}(\{\Omega_X^{\frac{n}{2}}\})\,\subset\, \Mf_X(n)
$$
which is isomorphic,
via $E\mapsto E_0=E\ot\Omega_X^{-\frac{1}{2}}$, to the moduli
of semistable $\text{SL}_n$--bundles. The dimension
of $\N_X(n)$ is $(g-1)(n^2-1)$. 

The subvariety
$$
\Theta \, :=\, \{V\in \Mf_X(n)\, \vert\, \Ho(X,\, V)\,
\not=\, 0\} 
$$
is a (reduced) divisor, the {\it generalized theta divisor},
that gives the ample generator of the Picard group
${\rm Pic}(\N_X(n))$ \cite{DN}. Note that
for any $E$ in $\Mf_X(n)$, we have $h^0(E)=h^1(E)$.
The condition $h^0(E)=h^1(E)=0$ also guarantees
that $E$ is semistable. Indeed, if a subbundle $F$ of $E$
contradicts the semistability condition of $E$, then the
Riemann--Roch theorem ensures that $h^0(F) >0$, thus
contradicting the condition that $h^0(E)=0$.
The smooth locus of the theta
divisor $\Theta$ is precisely the subvariety $\Theta^\circ$ of
vector bundles $E$
with ${\rm h}^0(E)={\rm h}^1(E)=1$. 

Let $\K_X(n)\subset\N_X(n)$ denote the subvariety consisting of vector
bundles, which are isomorphic to a direct sum of line bundles. Thus
for $n=2$, $\K_X(2)$ consists of vector bundles of the form $\Ll\oplus
\Ll^\vee\cong \Ll^{\vee}\oplus \Ll$, so that $\K_X(2)$ is isomorphic
to the Kummer variety $\K_X(2)=\Pg/\{\Ll\sim\Ll^{\vee}\}$.

\subsection{The Szeg\"o kernel.}
For $E\in \Mf_X(n)$, with $E_0=E\ot
\Omega_X^{-\frac{1}{2}}\in\Mf_X(n)_0$, denote by $\M(E)$ the sheaf
$$\M(E)\,=\,\M_1(E_0)\,\cong\, E\boxt E^{\vee}(\Delta)\, .$$
(By \remref{indep of theta} $\M(E)|_{n\Delta}$ is independent of
$\Ohalf$.) Let $\M(E)^\circ$ denote the subsheaf
$$\M(E)^{\circ}\,=\,\{s\in \M(E)\; : \;
s|_{\Delta}=\lambda\on{Id}_E\;(\la\in\C)\}.$$ When
$E\in\Mf_X(n)\sm\Theta$, there is a canonical kernel function
associated to $E$, the nonabelian Szeg\"o kernel of Fay
\cite{Fay nonabelian, Fay Szego}
(see also \cite{part 1}). In particular we will use the following
characterization of the Szeg\"o kernel:

\subsubsection{Proposition.}\label{kernels} (\cite{part 1})
\begin{enumerate}
\item If ${\rm h}^0(E)={\rm h}^1(E)=0$, then $\Ho(X\times
X,\M(E)^\circ)=\C\cdot\s_E$, where $\s_E$, the {\em Szeg\"o kernel} of
$E$, is the unique section with $\s_E|_{\Delta}=\on{Id}_E$.
\item Otherwise, the inclusion
$$
\Ho(X,E)\ot\Ho(X,E^\vee)\,\cong \,\Ho(\XX,E\boxt
E^{\vee})\,\hookrightarrow\, \Ho(\XX,\M(E)^\circ)
$$
is an isomorphism. In other words,
all global sections of $\M(E)^{\circ}$ vanish on $\Delta$. 
\end{enumerate}

\subsubsection{} Thus $\s_E|_{k\Delta}\in \Mop_k(E_0)(1)$ is a 
canonical (shifted) matrix oper on $E_0$ (\secref{mopers}).  The
proposition follows from Serre duality and the long exact sequence of
cohomologies of $E\boxtimes E^{\vee}$ with poles along the diagonal.

We would like to apply the determinant map to the Szeg\"o kernel:
$$
\det \s_E\,\in\, \Gamma(\M_n(\det E_0))\, .
$$
Restricting to $k\Delta$
defines a (shifted) $\text{GL}_k$--oper for the line bundle
$\det E_0$. (We will identify shifted opers with opers, using
\secref{shifted opers}.)

\subsubsection{Definition.} 
\begin{enumerate}
\item The Wirtinger oper associated to a bundle $E\in\Mf_X(n)\sm
\Theta$ is the $\text{GL}_n$--oper
$\det\s_E|_{(n+1)\Delta}\in \Gamma(\M_n(\det E_0)|_{(n+1)\Delta})$.
The resulting map
\begin{eqnarray*}
\W\,:\,\Mf_X(n)\sm\Theta &\longrightarrow& \Op_n\\
\W\,:\,\N_X(n)\sm\Theta&\longrightarrow &\Op_n^\circ
\end{eqnarray*} 
is the {\em Wirtinger map} (of rank $n$).
\item The Klein oper kernel associated to a bundle $E\in \Mf_X(n)\sm
\Theta$ is the kernel $\det \s_E\in \Ho(\XX,\M_n(\det E))$.
The resulting map
\begin{equation*}
\Kl\,:\,\N_X(n)\sm\Theta \,\longrightarrow\, \Kern_n
\end{equation*} 
is the {\em Klein map} (of rank $n$).

\end{enumerate} 

\subsubsection{} Note
that the dimensions of $\Mf_X(n)$ and $\Op_n$ agree, as do those of
$\N_X(n)$ and $\Op_n^\circ$. Thus if we knew $\W$ to be a finite map, it
would give a canonical system of \'etale coordinates on an open
subvariety of the moduli space. This leads us to conjecture:

\subsubsection{Conjecture.}
\begin{enumerate}
\item The Klein map is finite onto its image for all $X$.
\item The Wirtinger map is finite for generic $X$.
\end{enumerate}

\subsubsection{}
We will prove the conjecture in the case of torus bundles, i.e., along
$\K_X(n)\subset\N_X(n)$.
We first describe the Szeg\"o kernel and its determinant for
torus bundles. Suppose
$E\cong\Ll_1\oplus\Ll_2\oplus\cdots\oplus\Ll_n$. Then $E\in
\Mf_n(X)\sm \Theta$ if and only if each $\Ll_i\in \Pg\sm\Theta$.
Moreover in this case $\s_E=\s_{\Ll_1}\oplus\cdots\oplus \s_{\Ll_n}$,
and $\det \s_E =\bigotimes_{i=1}^n\s_{\Ll_i}$. If $E\in\N_X(n)$ then we
have in addition $\bigotimes_{i=1}^n\Ll_i=\Omega_X^{\frac{n}{2}}.$ For
example, if $n=2$, $E=\Ll\oplus\Ll^\vee$ and $\s_E=\s_\Ll
\s_{\Ll^{\vee}}=\s_\Ll\s_{\Ll}^t$.

Recall the {\em Petri map}
$$
\Ho(X,\Ll)\ot \Ho(X,\Ll^{\vee})\,\longrightarrow \,\Ho(X,\Omega)
$$
obtained by tensoring of sections \cite[p. 127]{ACGH}.
Under the K\"unneth isomorphism $$\Ho(X\times
X,\Ll\boxt\Ll^{\vee})\,=\,\Ho(X,\Ll)\ot \Ho(X,\Ll^{\vee})\, ,
$$
the Petri map is identified with the restriction to the diagonal
$$
\Ho(\XX,\Ll\boxt\Ll^{\vee})\,\longrightarrow\,
\Gamma(\Ll\boxt\Ll^{\vee}|_{\Delta})\,=\,\Ho(X,\Omega)\, .
$$
Thus injectivity of
the Petri map implies that global sections of $\Ll\boxt\Ll^{\vee}$
are determined by their restriction to the diagonal.
The curve $X$ is called {\em
Brill--Noether general} if the Petri map is injective for every line
bundle $\Ll$. By the Petri conjecture (Lazarsfeld's Theorem), this
condition is satisfied by a generic curve of genus $g$.

We then having the following result in the direction of the finiteness
conjecture: 

\subsubsection{Theorem.}\label{the theorem}
\begin{enumerate}
\item The Klein map for Kummers $\Kl:\K_X(n)\sm\Theta \to \Kern_n$ is
finite onto its image for all $X$.
\item The Wirtinger map for Kummers $\W:\K_X(n)\sm\Theta\to
\Op_n^\circ$ is finite onto its image for Brill--Noether general
$X$.
\end{enumerate}

\subsubsection{Proof of (1).}
Consider the subvariety of $(\Pg)^n$ of
line bundles $(\Ll_1,\cdots,\Ll_n)$ with $\bigotimes_{i=1}^n
\Ll_i\cong \Omega_X^{\frac{n}{2}}$. (We identify
this with $(\Pg)^{n-1}$ through the first
$n-1$ $\Ll_i$.) For $(1)$, it clearly suffices to
show that the map from $(\Pg)^{n-1}$ to $\Kern_n$ given by
$$
(\Ll_1,\dots,\Ll_{n-1})\,\longmapsto\, \bigotimes_{i=1}^n \s_{\Ll_i}
$$
is finite. To do
so we consider $\Kern_n$ as a subvariety of $\Pp\Ho(\XX,\M_n)$
(contained in the affine open of sections with nonzero trace on the
diagonal), and complete $\Kl$ to a morphism
$$
\Kl\,:\,(\Ptil)^{n-1}\,\longrightarrow\,
\Pp \Ho(\XX,\M_n)
$$
{}from a partial resolution of the singular locus of $\Theta$. 
Here $\Ptil\to \Pg$ is a projective morphism, which is an isomorphism
off the singular part of the theta divisor.
(In fact $\Ptil$ will be the union, 
for $X$ Brill--Noether general, of the projectivized conormal bundles
to the Brill--Noether loci $W^{g-1,i}\subset\Pg$.) Hence the extended map
is automatically proper, and a closer examination shows it remains proper
when restricted to $(\Pg\sm\Theta)^{n-1}$, and hence finite.

We construct $\Ptil$ as the moduli of pairs $(\Ll,s)$ consisting of a
line bundle $\Ll\in\Pg$ and a nonzero section $s$ of $\M(\Ll)$, up to scalar
(i.e., a divisor in the
complete linear series $|\M(\Ll)|$ on $\XX$.) 
This is a projective variety mapping to $\Pg$, with the fibers the
projective spaces $\Pp \Ho(\XX,\Ll\boxt\Ll^{\vee}(\Delta))$.  The
construction of this projective variety follows from that of the
Hilbert scheme of divisors, of the same degree as $\M(\Ll)$, on the
surface $X\times X$. This Hilbert scheme fibers over the Picard group
of $\XX$, and we pull it back to $\Pg$ over the morphism
$\Pg\to\on{Pic}(\XX)$ sending $\Ll$ to $\M(\Ll)$.

It follows from \propref{kernels} that over $\Pg\sm\Theta$ the
projection $\Ptil\to \Pg$ is an isomorphism, since the Szeg\"o kernel
is the unique section of $\M(\Ll)$ up to scalars. In fact, the morphism
remains an isomorphism on the smooth locus of $\Theta$, since for
$h^0(\Ll)=1$ we have $h^0(\M(\Ll))=h^0(\Ll)h^0(\Ll^\vee)=1.$ Since by
\propref{kernels} every section of $\M(\Ll)$ for $\Ll\in\Theta$
defines a section of $\Ll$ and one of $\Ll^{\vee}$, it follows that
the inverse image in $\Ptil$ over $\Theta$ (for the
projection of $\Ptil$ to $\Pg$) is given by
$$
\Ptil|_{\Theta}\,\cong\, \on{Sym}^{g-1}X\times_{\Pg}
i^*\on{Sym}^{g-1}X\, ,
$$
where $i:\Ll\to\Ll^{\vee}$ -- in other words, the inverse image is the
space of pairs of divisors for $\Ll$ and $\Ll^{\vee}$. (Thus $\Ptil$
restricts, for $X$ Brill--Noether general, to the union of blowups of
the Brill--Noether loci in $\Pg$.)

We now extend the morphism $\Kl$ from $(\Pg)^{n-1}$ to $\Ptil_n$, the
inverse image of $(\Pg)^{n-1}\subset (\Pg)^n$ in $(\Ptil)^n$,
i.e., $\Ptil_n$ parametrizes $(\Ll_1,s_1;\dots;\Ll_n,s_n)$ where the
$\Ll_i$ add up to $\Omega_X^{\frac{n}{2}}$. To such a tuple we assign
the line $[\bigotimes_{i=1}^n s_i]$ in
$$
\bigotimes_{i=1}^n (\pi_{\XX})_*\M(\Ll_i)\,=\,
(\pi_{\XX})_*\M(\Omega_X^{\frac{n}{2}})\, ,
$$ 
where $s_i$ are the tautological sections of $\Ll_i$ given
by the $i$th point in $\Ptil$ (taken up to scalar). The
right hand side is the vector space $\Ho(\XX,\M_n)$, independently of
the $\Ll_i$, so we have constructed the desired extension
$$
\Kl\,:\,\Ptil_n\,\longrightarrow \,\Pp\Ho(\XX,\M_n)\, .
$$

The completed morphism $\Kl$ is a morphism of projective varieties,
hence proper. We claim its restriction to $(\Pg\sm\Theta)^{n-1}$ is
also proper. Let
$$
\Pp\Ho(\XX,\,\M_n)\,\subset\, \Pp\Ho(\XX,\,\M_n)
$$
denote the hyperplane of sections vanishing on the diagonal. By
\propref{kernels}, for $\Ll\in \Theta$, all sections of $\M(\Ll)$
automatically vanish on the diagonal, while for $\Ll\in\Pg\sm\Theta$
all nonzero sections give nonzero constant functions on the diagonal.
Hence the preimage of the complement of this hyperplane is precisely
$\Pg\sm \Theta$. We obtain that the morphism $\Kl$ from the affine
variety $(\Pg\sm\Theta)^{n-1}$ is proper, hence finite.

\subsubsection{Proof of (2).}
We embed the affine space $\Op_n^\circ$ in the projective space
$$\ol{\Op}_n^\circ \,=\, \Pp\Gamma(\M_n|_{(n+1)\Delta})\, .$$
Thus $\W$ gives rise to a map
$$
\W\,:\, (\Pg\sm\Theta)^{n-1}\,\longrightarrow\, 
\ol{\Op}_n^\circ\, .
$$
In order to
prove finiteness of $\W$, we would like to extend it
to $\Ptil$, whenever possible.

Let $\Ll\in\Theta$. Then by \propref{kernels}, global sections of
$\Ll\boxt\Ll^{\vee}(\Delta)$ vanish on $\Delta$. If the Petri map of
$\Ll$ is injective, however, such sections are determined by their
restriction to $2\Delta$. So we take $X$ to be Brill--Noether general.
It follows that for a collection of nonzero sections $s_i\in
\Ho(\XX,\M(\Ll_i))$, the restriction $(\bigotimes_{i=1}^n
s_i)|_{(n+1)\Delta}$ is also nonzero. Thus
the $s_i$ define a point in $\ol{\Op}_n^\circ$, and we have completed
$\W$ to a map
$$
\W\,:\,(\Ptil)^{n-1}\,\longrightarrow\, \ol{\Op}_n^\circ\, .
$$
Again the
inverse image of the hyperplane of sections vanishing on the diagonal
is precisely the inverse image of the theta divisor, so the map
remains proper off $\Theta$, implying finiteness as before.

\section{Relations with theta functions.}\label{thetas}

\subsection{The theta linear series.}\label{2theta}
The Klein and Wirtinger maps have natural interpretations as quotients
of the theta linear series on $\Mf_X(n)$ and $\N_X(n)$. 
For $E\in\Mf_X(n)$, consider the sequence of maps
$$
((n+1)\Delta)\,\hookrightarrow\, X\times
X\, \stackrel{\delta}{\longrightarrow}\, \Jac
\,\stackrel{\tau_E}{\longrightarrow}\, \Mf_X(n)\, .
$$ 
Here
$$\tau_E\,:\,\Jac \,:=\,\text{Pic}^0(X)\, \longrightarrow\,
\Mf_X(n),\hskip.3in
\tau_E(\Ll)\,=\,E\ot\Ll
$$
is the translation map, $\delta(x,y)=y-x$ 
and the composition $\tau_E\circ\delta$ is the difference map
$$
\delta_E\,:\,\XX\,\longrightarrow \,\Mf_X(n),\hskip.3in
\delta_E(x,y)\,=\,E(y-x)\, .
$$

It is well--known that for $E\in\N_X(n)$, the
pullback of nonabelian theta functions
$$\tau_E^*[\Oo_{\Mf_X(n)}(\Theta)]\,=\,\Oo_{\Jac}(n\Theta)$$
are weight $n$ abelian theta functions. Moreover the resulting map
$$
\tau^*\,:\,\N_X(n)\,\longrightarrow\, \Pp\Ho(\Jac,\Oo_{\Jac}(n\Theta))
$$
is an embedding (see \cite{Beauville}). (Note that we have fixed a
theta characteristic $\Ohalf$, which allows us to principally polarize
the Jacobian and pass from line bundles $\Ll$ of degree $n(g-1)$ to
$\Ll_0$ of degree $0$.) Pulling back further to $X\times X$ or
$(n+1)\Delta$, we obtain sections of the pullback
$\delta_E^*\Oo_{\Mf_X(n)}(\Theta)=\M_n\otimes\Theta|_E,$ the tensor of
the line bundle $\M_n$ by the complex line $\Theta|_E$, the fiber of
$\Theta$. (See e.g. \cite{part 1}.)

It follows that we have a sequence of pullback maps
$$
\Ho(\Jac,\,\Oo_{\Jac}(n\Theta))\,\longrightarrow\, \Ho(\XX,\,\M_n)
\,\longrightarrow\, \Gamma(\M_n|_{(n+1)\Delta})\, ,
$$
and consequently
rational maps on the corresponding projective spaces. Composing these
with $\tau^*$ we obtain rational maps from $\N_X(n)$ (if the image
of $\tau^*$ is not contained in the kernels of the projections).

We
will use the following description of the Szeg\"o kernel:

\subsubsection{Theorem.}\label{kernel and theta} (\cite{part 1}, see also \cite{GP,Po})
 $\on{det}\s_E\,=\,\delta_E^*\theta/\theta(E)$.

\subsubsection{Corollary.} The Klein and Wirtinger maps
\begin{eqnarray*}
\Kl\,:\,\N_X(n)\sm\Theta&\longrightarrow& \Pp\Ho(\XX,\M_n)\\
\W\,:\,\N_X(n)\sm\Theta&\longrightarrow &
\Pp\Gamma(\M_n|_{(n+1)\Delta})
\end{eqnarray*}
are equal to the composition of the theta linear series $\tau^*$ with
the restrictions to $\XX$ and $(n+1)\Delta$, respectively.

\subsection{The linear series $|2\Theta|$.} 
Let us consider the case $n=2$. (Our reference for $2\Theta$ functions
is \cite{Do}.) 
The map $\tau^*:\N_X(n)\to\Pp\Ho(\Jac,\Oo(2\Theta))$ restricts on the
Kummer variety $\Jac\twoheadrightarrow \K_X(2)\subset \N_X(2)$ to the map
$$
\Jac\,\ni\, e\,\longmapsto \,\Theta_{e}+\Theta_{-e}
$$
(where $\Theta_e$ denotes the translate of $\Theta$ by $e$).
The Riemann quadratic identity and
Kummer identification theorem provide a natural isomorphism between
this map and the $2\Theta$ linear series
$$
|2\Theta|_*\,:\,\Jac\,\longrightarrow\, \Pp\Ho(\Jac,\Oo(2\Theta))^*
$$
which naturally maps to the dual projective space.

By the symmetry properties of $2\Theta$ it follows that the image of
$\Ho(\Jac,\Oo(2\Theta))$ in $\Ho(\XX,\M_2)$ consists of {\em symmetric}
bidifferentials. In fact 
there is a short exact sequence
\begin{equation*}
0\,\longrightarrow\,\Gamma_{00}\,\longrightarrow\,
\Ho(\Jac,\,\Oo(2\Theta))\,\stackrel{\delta^*}{\longrightarrow}\,
\Ho(\XX,\,\Omega\boxt\Omega(2\Delta))^{\text{sym}}
\,\longrightarrow\, 0\, ,
\end{equation*}
where the kernel $\Gamma_{00}$ can be characterized as the subspace of
$2\Theta$--functions vanishing to fourth order at $0$. The right hand
side is a vector space of dimension $\binom{g}{2}+1$. Its
projective
space $\ol{\Kern_2^{sym}}\cong \Pp^{\binom{g}{2}}$ contains as an
affine open the space $\Kern_2^{sym}$ of projective kernels. This
vector space has a further quotient
$\Gamma(\Omega_X\boxt\Omega_X(2\Delta)|_{3\Delta})^{\text{sym}}$,
obtained by
restricting kernels to $3\Delta$. Its projective space
$\ol{\Proj}\cong \Pp^{3g-3}$ contains as an affine open the space
$\Proj$ of projective structures. Note that the image of
$\Kl:\N_X(2)\sm\Theta\to\Kern_2$ lies in $\Kern_2^{sym}$, while $\W$
defines a map $\W:\N_X(2)\sm \Theta\to \Proj$. We may thus reinterpret
the finiteness theorem as follows:

\subsubsection{Corollary.} 
\begin{enumerate}
\item The rational map 
$\Jac\to \Pp^{\binom{g}{2}}$ defined by the composition of
$|2\Theta|_*$ with projection by $\Gamma_{00}$ is a finite
morphism on $\Jac\sm\Theta$.
\item For $X$ generic, the further projection $\Jac\to\Pp^{3g-3}$ 
remains finite on $\Jac\sm\Theta$.
\end{enumerate} 

\subsubsection{Formulas.}
The explicit description of the Szeg\"o kernel for line bundles is
\begin{equation}\label{formula for Szego}
\s_\Ll(x,y)\,=\,\frac{\theta(y-x+\Ll_0)}{\theta(\Ll_0)E(x,y)}\, ,
\end{equation}
where $E(x,y)$ is the prime form (this is the rank one case of
\thmref{kernel and theta}.)
Thus the Klein map on the Kummer $\K_X(2)$ becomes
\begin{equation}\label{Sug in rank one}
\Kl(\Ll\oplus \Ll^{\vee})\,=\, \s_{\Ll}\s_{\Ll^{\vee}}\,=\,
\frac{\theta(y-x+\Ll_0)\theta(y-x-\Ll_0)}
{\theta(\Ll_0)^2 E(x,y)^2}\, .
\end{equation}

The relation to $2\Theta$ is easily seen explicitly.
Let
$$
\arrtheta\,:\,\C^g\,\longrightarrow\, \Ho(\Jac,\Oo(2\Theta))^*
$$
defined by
$$
\arrtheta(e)\,=\, \sum_{\alpha,\beta\in\Jac[2]}
\theta_2\left[\begin{array}{c} \alpha\\ \beta\end{array}\right](e)
$$
be the generating vector of the second order theta functions with
characteristics. By Riemann's quadratic identity (\cite{Tata I}), we
may rewrite the expression $\Kl(\Ll\oplus \Ll^{\vee})$ of \eqref{Sug
in rank one} as follows:
\begin{equation}\label{using Gauss}
\frac{\theta(y-x+\Ll_0)\theta(y-x-\Ll_0)}
{\theta(\Ll_0)^2 E(x,y)^2}\,=\,
\frac{\arrtheta(y-x)\cdot\arrtheta(\Ll_0)}{\theta(\Ll_0)^2E(x,y)^2}\, .
\end{equation}

\subsection{The Gauss map.}
Let $\Theta^{\circ}\subset\Pg$ denote the smooth part of the theta
divisor. The Gauss map for the theta divisor sends
$$
{\mf G}\,:\,\Theta^{\circ}\,\longrightarrow\, \Pp \Ho(X,\,\Omega)\, .
$$
Since $\Ho(X,\,\Ll)\,=\,\C l$
is one dimensional for $\Ll\in\Theta^{\circ}$, the
Petri map for $\Ll$,
$$
\Ll\,\longmapsto\, l\ot l^{\vee}\, ,
$$
also defines a line
in $\Ho(X,\Omega)$, which is known to agree with the Gauss line for
$\Ll$. On the other hand the extension of $\W$ to
$\Theta^{\circ}\subset \Ptil$ sends
$$
\Ll\,\longmapsto\, (l\boxt l^{\vee})\ot
(l^{\vee}\boxt l)|_{\Delta}\,=\,(l\ot l^{\vee})^{\ot 2}\, ,
$$
which defines a line in $\Pp\Ho(X,\Omega)\subset\ol{\Proj}$.
Thus the tensor square of the Gauss map agrees with the morphism $\W$:

\subsubsection{Corollary.} For a Brill--Noether general curve, 
the square of the Gauss map
$${\mf G}^{\ot 2}\,:\,\Theta^{\circ}\,\longrightarrow\,
\Pp\Ho(X,\Omega^{\ot 2})
$$ extends to a finite morphism
$$\W\,:\,\Pg\sm\Theta^{sing}\,\longrightarrow \,\ol{\Proj}\, .$$

\subsubsection{Remark.} It is interesting to note that this 
relation of the Klein map to the theta divisor fails completely in
higher rank. Namely, for $E\in\Theta^{\circ}$ we still have
$\Ho(X,E)=\C s$. It follows that the Higgs field
$$
s\ot s^\vee\,\in\,
\on{End}E\ot\Omega\,=\, (E\boxt E^{\vee})|_{\Delta}
$$
is {\em nilpotent}. In
fact as $E$ varies over $\Theta^{\circ}$ we obtain this way an
irreducible component of the global nilpotent cone in the moduli of
Higgs bundles. Thus the ``Hitchin--Gauss'' map, applying
characteristic polynomials to this canonical line of Higgs bundles
along $\Theta^\circ$, is identically zero. In particular the
determinant $\det s\boxt s^{\vee}=0$ vanishes identically on $\XX$, so
we cannot use this to extend the Klein map across the theta divisor.

\subsection{Logarithmic derivatives of theta.}
In \cite{Tata I}, Mumford cites three general techniques for
constructing meromorphic functions on Jacobians out of theta
functions, of which the third is that of taking second logarithmic
derivatives. Namely, there is a collection of $\binom{g}{2}$
meromorphic functions
$$\frac{\partial^2\log
\theta}{\partial z_i \partial z_j}$$ on the Jacobian -- or more
invariantly, a rational map 
\begin{equation}\label{ddlog theta}
e\,\longmapsto\, \sum_1^g \frac{\partial^2\log
\theta}{\partial z_i \partial z_j}(e) \omega_i(x)\omega_j(y)\, .
\end{equation}
{}from the Jacobian to holomorphic symmetric bidifferentials on $X$,
$\Ho(\XX,\Omega\boxt\Omega)^{\text{sym}}$.

By translating these holomorphic bidifferentials by the Bergman kernel
$\omega_B$ (\secref{Bergman}), we obtain the {\em Klein projective kernels} 
$\omega_e\in\Kern_2^{sym}$ (\cite{Tyurin}):
\begin{equation}\label{Klein bidiff}
\omega_e\,=\,\omega_B(x,y)+\sum_1^g \frac{\partial^2\log
\theta}{\partial z_i \partial z_j}(e) \omega_i(x)\omega_j(y)\, .
\end{equation}
Here the point $e\in\Jac\sm\Theta$. Classically $e$ is taken to be a
two--torsion point, so that $\omega_e$ is written in terms of theta
functions with characteristics. The corresponding projective
connections $\omega_e|_{3\Delta}\in\Proj$ are known (\cite{Tyurin}) as
the Wirtinger connections.

The relation of these classical kernels with our Klein and Wirtinger
maps is provided by the ``second corollary to the trisecant
identity'' of J. Fay (\cite{Fay}, Corollary 2.12; also \cite{Tata
II}):
\begin{equation*}
\Kl(\Ll)=\frac{\theta(y-x+\Ll_0)\theta(y-x-\Ll_0)} {\theta(\Ll_0)^2
E(x,y)^2}=\omega_B(x,y)+\sum_1^g \frac{\partial^2\log \theta}{\partial
z_i \partial z_j}(\Ll_0) \omega_i(x)\omega_j(y)=\omega_{\Ll_0}\, .
\end{equation*}

\subsubsection{Corollary.}\label{log theta} The second logarithmic
derivatives of $\theta$ provide a finite parametrization of the
complement of the theta divisor in the Jacobian in affine space of
dimension $\binom{g}{2}$. Namely, the holomorphic map 
$$
\Jac\sm\Theta\,\longrightarrow \,\Ho(\XX,\,\Omega\boxt
\Omega)^{\text{sym}}
$$
of \eqref{ddlog theta} is finite onto its image.

\subsubsection{Remark.} 
It also follows from \corref{log theta} that the second logarithmic
derivative map is generically finite for generic abelian varieties,
since it is finite on the Jacobian locus.

\medskip
\noindent {\bf Acknowledgments:}\, We would like to thank Ron Donagi
and Matthew Emerton for useful discussions. We are especially grateful
to Mohan Ramachandran for suggesting the use of properness to
establish finiteness.


\end{document}